\documentclass[11pt]{amsart}
\usepackage{amsthm,amsmath,amssymb,amscd,epsfig}

{\end{list}}
{
   \newtheorem{theorem}{Theorem}[section]
   
   \newtheorem{lemma}[theorem]{Lemma}

}
{\theoremstyle{definition}

}

\newcommand{\RR}{{\mathbb{R}}}

\newcommand{\QQ}{{\mathbb{Q}}}

\newcommand{\ZZ}{{\mathbb{Z}}}

\newcommand{\cF}{{\mathcal F}}

\newcommand{\cO}{{\mathcal O}}

\newcommand{\Relint}{\operatorname{Rel.Int}}

\newcommand{\Coker}{\operatorname{Coker}}

\newcommand{\Span}{\operatorname{Span}}

\newcommand{\Hom}{{\operatorname{Hom}}}

\newcommand{\setmin}{{\smallsetminus}}

\begin{document}
\title{On toric vector bundles}

\author{Saman Gharib}
\address{Department of Mathematics\\ University of British Columbia \\
  1984 Mathematics Road\\
Vancouver, B.C. Canada V6T 1Z2}
\email{samangharib50233@gmail.com}

\author{Kalle Karu}
\address{Department of Mathematics\\ University of British Columbia \\
  1984 Mathematics Road\\
Vancouver, B.C. Canada V6T 1Z2}
\email{karu@math.ubc.ca}

\thanks{This work was partially supported by NSERC Discovery and USRA grants.}

\begin{abstract} 
Following Sam Payne's work, we study the existence problem of nontrivial vector bundles on toric varieties. The first result we prove is that every complete fan admits a nontrivial conewise linear multivalued function. Such functions could potentially be the Chern classes of toric vector bundles. Then we use the results of  Corti\~{n}as, Haesemeyer, Walker and Weibel to show that the (non-equivariant) Grothendieck group of the toric 3-fold studied by Payne is large, so the variety has a nontrivial vector bundle. Using the same computation, we show that every toric 3-fold X either has a nontrivial line bundle, or there is a finite surjective toric morphism from Y to X, such that Y has a large Grothendieck group.  

\vspace{.3in}

{\bf CORRECTION.} One of the main results in this paper contains a fatal error. We cannot conclude the existence of nontrivial vector bundles on $X$ from the nontriviality of its $K$-group. The $K$-group that is computed here is the Grothendieck group of perfect complexes and not vector bundles. Since the varieties are not quasi-projective, existence of nontrivial perfect complexes says nothing about the existence of nontrivial vector bundles. We thank Sam Payne for drawing our attention to the error and Christian Haesemeyer for explanations about the K-theory.

\end{abstract}

\maketitle

\section{Introduction}

The purpose of this article is to study the existence of nontrivial vector bundles on toric varieties. Sam Payne in \cite{Payne} raised the question whether any complete toric variety has a nontrivial vector bundle. He gave an example of a $3$-dimensional toric variety that does not admit a non-trivial vector bundle of rank at most $3$. His proof was based on studying the equivariant Chern classes of torus-equivariant bundles (toric bundles for short). Since a toric vector bundle splits on a torus invariant affine open subset, the equivariant total Chern classes also split on such open sets. Payne proved that the particular toric $3$-fold does not admit any nontrivial equivariant Chern classes of degree $3$ or less, and toric vector bundles of rank $3$ or less with trivial equivariant total Chern class are all trivial.

We will first show that every complete toric variety has a nontrivial equivariant Chern class of sufficiently high degree. (By a Chern class we mean a cohomology class that splits on every torus invariant affine open, hence is a candidate for the total equivariant Chern class of a toric vector bundle.) The classes that we construct are such that they deform together with the fan.

To study the existence of vector bundles, we look at the (non-equivariant) $K$-theory of toric varieties. The equivariant $K$-theory of smooth toric varieties is well understood (\cite{VV}). Singular toric varieties, on the other hand, can have very complicated and large $K$-groups. Joseph Gubeladze in \cite{Gubeladze} gave an example of a toric variety $X$ such that $K_0(X)$ has an uncountable rank. Another such example was given by Corti\~{n}as, Haesemeyer, Walker and Weibel in \cite{CHWW}. We will use the results of \cite{CHWW} to show that the  toric $3$-fold studied by Payne is of the same type: its $K$-group has an uncountable rank. This in particular shows that the $3$-fold has nontrivial (non-equivariant) vector bundles of high enough rank. However, we were unable to give an explicit construction of a nontrivial bundle on this toric variety.

Using the same K-group computation, we also show that every complete toric $3$-fold $X$ either has a nontrivial line bundle, or there exists a finite  surjective toric morphism $X'\to X$ such that $K_0(X')$ has uncountable rank. The fact that a finite cover has big $K$-group is not surprising, because the cover is usually very singular. However, the fact that the cover has nontrivial $K$-group is not obvious to us.

We do not settle the question of whether there exists a complete toric variety with no nontrivial vector bundles. The results proved in the paper point to the opposite, that every complete toric variety admits a nontrivial vector bundle.

\section{Equivariant Chern classes of toric vector bundles.}

We refer to the book by Fulton \cite{Fulton} for background on toric varieties. We work over an algebraically closed field $k$ of characteristic zero. Let $N$ be a lattice and $M$ its dual lattice. A toric variety is determined by a pair $(\Delta, N)$, where $\Delta$ is a fan in the lattice $N$. The set of one-dimensional cones of $\Delta$ is denoted $\Delta(1)$. If $\rho\in\Delta(1)$, then $n_\rho$ is the first nonzero lattice point on $\rho$. The affine toric variety associated to a cone $\sigma$ in $N$ is $U_\sigma$, the big torus is $T$, and the divisor corresponding to $\rho\in\Delta(1)$ is denoted $V_\rho$.  

We refer to \cite{Payne} for notation and background material on toric vector bundles. A toric vector bundle is a torus equivariant vector bundle on a toric variety. Such bundles were classified by Klyachko \cite{Klyachko}.

Let $X_\Delta$ be a toric variety corresponding to a fan $\Delta$. A toric line bundle on $X_\Delta$ is determined by its equivariant first Chern class, which is a conewise linear integral function on the support of the fan $\Delta$. A toric vector bundle splits into a direct sum of line bundles on every affine open $U_\sigma \subset X_\Delta$ for $\sigma\in \Delta$. This splitting gives an $r$-element  multiset $L_\sigma$ of integral linear functions on every cone $\sigma$. When $\tau$ is a face of $\sigma$, then the multiset $L_\sigma$ restricts to the multiset $L_\tau$. The collection $\{L_\sigma\}_{\sigma\in\Delta}$ is called a multivalued integral conewise linear function on $\Delta$. Its degree is $r=|L_\sigma|$. A multivalued conewise linear function is called trivial if there exists a multiset $L$ of global linear functions, such that $L_\sigma$ is the restriction of $L$ to $\sigma$ for any cone $\sigma\in\Delta$.

Given a symmetric polynomial $s(x_1,\ldots,x_r)$ in $r$ variables and a multivalued conewise linear function $\{L_\sigma\}_{\sigma\in\Delta}$ of degree $r$, one can produce a conewise polynomial function by substituting on each cone $\sigma$ the functions in $L_\sigma$ into $s$. When $s_i$ is the $i$-th elementary symmetric function and $\{L_\sigma\}$ the multivalued conewise linear function associated to a toric vector bundle $E$, then the conewise polynomial function produced this way is the $i$-th equivariant Chern class of $E$. In other words, the multivalued function encodes all the equivariant Chern class data of $E$.

Ignoring the existence of toric vector bundles for a moment, we can ask whether every complete fan  $\Delta$ admits a nontrivial multivalued conewise linear  integral function. This is indeed true.

\begin{theorem}\label{eq:main} Let $\Delta$ be a fan (rational, polyhedral) in the lattice $N$ having more than one maximal cone with dimension equal to the rank of $N$. Then there exists a nontrivial multivalued conewise linear integral function on $\Delta$.
\end{theorem}

\begin{lemma}\label{eq:main-lem}
For every rational cone $\sigma$ there exist two different multi-valued linear integral functions on it such that the restriction of these two functions on the facets of $\sigma$ are equal as multi-valued linear functions.
\end{lemma}

{\em Proof of the lemma.} Let $\tau_1,\ldots,\tau_k$ be the facets of $\sigma$.  Pick nonzero elements $L_i\in M\cap \sigma^\vee\cap \tau_i^\perp$. For $j=1,2$, consider the $2^{k-1}$-element multisets
\[ A_j = \{ \varepsilon_1 L_1 + \cdots +\varepsilon_k L_k| \varepsilon_i\in\{0,1\}, \sum_i\varepsilon_i= j \pmod{2}\}.\]
Then the two multisets $A_1, A_2$ restrict to the same multiset on any $\tau_i$ (to get a bijection between the two restrictions to $\tau_i$, just change the coefficient $\varepsilon_i$).  The two multisets $A_1, A_2$ are not equal; for example, $0$ lies in $A_2$, but not in $A_1$. \qed

{\em Proof of the theorem.}  Pick an arbitrary  cone $\sigma\in\Delta$ of maximal dimension. By the lemma, there are two distinct multivalued linear functions $A_1, A_2$ on $\sigma$ that restrict to the same multivalued function on every face $\tau$ of $\sigma$. To define the conewise linear function $L$ on $\Delta$, set $L$ equal to $A_1$ on $\sigma$ and to $A_2$ outside of $\sigma$. This function is nontrivial, i.e., not a multiset of global linear functions, because there exists another  cone $\sigma'$ of maximal dimension on which the two multisets $A_1, A_2$ are not equal. \qed 

Note that if we continuously deform the fan $\Delta$, we can also continuously deform the functions $L_1,\ldots,L_k$ in the lemma (of course, the functions will not be integral in general). This implies that the nontrivial conewise linear multivalued function on $\Delta$ also deforms with the fan. 

The next step in constructing an equivariant  vector bundle would be to fix a equivariant Chern class as in the theorem and look for a vector bundle corresponding to it. Such vector bundles may not exist. For example, if $\Delta$ is the fan over the faces of the regular cube with vertices $(\pm1,\pm1,\pm1)$ and $L$ the multivalued conewise linear function on $\Delta$ constructed in the theorem (for any choice of $\sigma$ and linear functions $L_1,L_2,L_3, L_4$), then one can check from Klyachko's description of equivariant bundles that there is no vector bundle with this Chern class. We will instead turn to $K$-theory to find nontrivial bundles.

\section{The Grothendieck's $K$-group.}
 
We use the results of Corti\~{n}as, Haesemeyer, Walker and Weibel \cite{CHWW} to show that the example $3$-fold studied by Payne \cite{Payne} has a large Grothendieck group, hence it must have nontrivial vector bundles. For general complete toric $3$-folds, we show that it either has a nontrivial line bundle or else we can construct a finite surjective cover with large Grothendieck group. 

The groups $K_0(X)$ studied here are non-equivariant. Since the equivariant $K$-group $K^T_0(X)$ may not surject onto $K_0(X)$,  the fact that $K_0(X)$ is nontrivial does not imply that there exists a nontrivial equivariant vector bundle on $X$.

Let $X=X_\Delta$ be a toric $3$-fold. It is shown in \cite{CHWW} (among several more general results) that the rational $K$-group $K_0(X)\otimes \QQ$ contains as a direct summand the cohomology $H^1(X,\cF)$ of a coherent sheaf $\cF$ on $X$. This sheaf is the cokernel of a morphism $\Omega^1_X \to \tilde\Omega^1_X$, where $\Omega^1_X$ is the sheaf of regular differential $1$-forms on $X$,  $\tilde\Omega_X^1$ is the Danilov's sheaf of $1$-forms (described below). Since the sheaves are defined over the ground field $k$, the cohomology group is a $k$-vector space. Thus, if $k$ is an extension of $\QQ$ of uncountable degree and the cohomology group is nonzero, then the rank of $K_0(X)$ must also be uncountable. We will check below that the $3$-fold in \cite{Payne} does indeed have a nonzero cohomology group. For certain more general $3$-folds we find finite covers where the cohomology group is nonzero.  
 
 Let us first recall the definitions of the sheaf $\tilde\Omega^1_X$ and the morphism $\Omega^1_X \to \tilde\Omega^1_X$ following \cite{CHWW}.
The sheaf $\tilde\Omega^1_X$ is defined by the exact sequence
\[ 0\to \tilde\Omega^1_X \to \cO_X\otimes M \stackrel{\delta}{\to} \oplus_{\rho\in\Delta(1)} \cO_{V_\rho},\]
where $\delta$ maps $f\otimes q \in \cO_X\otimes M$ to $f|_{V_\rho}\cdot \langle n_\rho, q\rangle$. The sheaf $\cO_X\otimes M$ can be identified with the sheaf $\Omega^1_X(\log D)$ of differentials with log poles along $D=X\setmin T$. The identification takes $f\otimes q \in \cO_X\otimes M$ to the form $f\cdot \frac{d\chi^q}{\chi^q}$. 
 
 The sheaf $\tilde\Omega^1_X$ is $T$-equivariant, hence its sections on open $T$-invariant sets $U_\sigma$ are $M$-graded.  If $m\in \sigma^\vee\cap M, q\in M$, then $\chi^m \otimes q$ is homogeneous of degree $m$. The sections $\tilde\Omega^1_X(U_\sigma)$ of degree $m$ can be described as follows.  If $m\notin \sigma^\vee$, then  $\tilde\Omega^1_X(U_\sigma)_m=0$. Otherwise, let $\sigma^\vee_m$ be the smallest face of $\sigma^\vee$ containing $m$. Then
 \[ \tilde\Omega^1_X(U_\sigma)_m = k\chi^m \otimes (\Span(\sigma^\vee_m)\cap M).\]
 
 The morphism $\Omega^1_X \to \tilde\Omega^1_X$ is defined by sending a form $\chi^p d\chi^q \in \Omega^1_X(U_\sigma)$ to 
 \[ \chi^{p+q} \frac{d\chi^q}{\chi^q} = \chi^{p+q}\otimes q \in \cO_X(U_\sigma)\otimes M.\]
 This image lies in $\tilde\Omega^1_X$. 
 
 The sheaf $\Omega^1_X$ is again $T$-equivariant. The section $\chi^p d\chi^q \in \Omega^1_X(U_\sigma)$ is homogeneous of degree $p+q$. It follows that the morphism $\Omega^1_X \to \tilde\Omega^1_X$ is also $T$-equivariant. Let us identify the image of the map $\Omega^1_X(U_\sigma)_m \to \tilde\Omega^1_X(U_\sigma)_m$. If $p,q\in \sigma^\vee\cap M$, $p+q=m$, then $\chi^p d\chi^q$ has degree $m$ and maps to $\chi^{p+q} \otimes q$. Thus, the image is
 \[ k \chi^m\otimes (\Span\{q | q\in \sigma^\vee\cap M, m-q\in\sigma^\vee\cap M\}\cap M) \subset k\chi^m \otimes (\Span(\sigma^\vee_m)\cap M).\]
 A good way to visualize this image is to note that 
 \[ \{q | q\in \sigma^\vee\cap M, m-q\in\sigma^\vee\cap M\} = \sigma^\vee \cap (-\sigma^\vee+m)\cap M.\]
 A particular case that we need below is when $\Span(\sigma^\vee_m)$ has dimension $1$. In that case 
 $\tilde\Omega^1_X(U_\sigma)_m$ has dimension $1$ and the image of $\Omega^1_X(U_\sigma)_m$ is nonzero (because we can take $q=m$), hence $\Omega^1_X(U_\sigma)_m \to \tilde\Omega^1_X(U_\sigma)_m$ is surjective.
 
 \subsection{Payne's toric $3$-fold}
 We now turn to the $3$-fold in \cite{Payne}. Let $N=\ZZ^3$ and start with the fan over the faces of the regular cube in $\RR^3$ with vertices $(\pm1, \pm1,\pm1)$. To get $\Delta$, we deform this fan by  replacing $(1,-1,1)$ with $(1,-1,2)$ and $(1,1,1)$ with $(1,2,3)$. 
 
 \begin{lemma} With $\Delta$ as above, the group $K_0(X(\Delta))$ has uncountable rank.
 \end{lemma}
 
 \begin{proof}
 Let $\tau\in \Delta$ be the $2$-dimensional cone generated by $(1,-1,-1)$ and $(1,-1,2)$, and let $\sigma_1$, $\sigma_2$ be the two $3$-dimensional cones having $\tau$ as a face.   
 
First let $X= U_{\sigma_1}\cup U_{\sigma_2}$. Denote $\cF = \Coker(\Omega^1_X \to \tilde\Omega^1_X)$. Mayer-Vietoris sequence for the cover 
 $X= U_{\sigma_1}\cup U_{\sigma_2}$ is:
 \[ \ldots \to H^0(U_{\sigma_1}, \cF)\oplus H^0(U_{\sigma_2},\cF) \to H^0(U_\tau, \cF)\to H^1(X,\cF) \to  H^1(U_{\sigma_1}, \cF)\oplus H^1(U_{\sigma_2},\cF) \to \ldots.\]
 Since $U_{\sigma_i}$ are affine and $\cF$ is coherent, the cohomology groups $H^1(U_{\sigma_i}, \cF)$ vanish for $i=1,2$. 
 
 Let us fix $m=(1,-1,0)\in \ZZ^3=M$. We claim that $H^0(U_{\sigma_i},\cF)_m = 0$ for $i=1,2$, but 
    $H^0(U_\tau, \cF)_m\neq 0$. This implies that $H^1(X,\cF) \neq 0$.
    
Let $\sigma_1$ be the cone generated by $(1,-1,-1),(1,-1,2), (1,1,-1), (1,2,3)$. Since $\langle m, (1,2,3)\rangle <0$, $m\notin\sigma^\vee$, hence $H^0(U_{\sigma_i},\cF)_m = 0$. The other cone $\sigma_2$ is generated by $(1,-1,-1),(1,-1,2), (-1,-1,-1), (-1,-1,1)$. Now $m\in\sigma_2^\vee$, but it vanishes on the face of $\sigma$ generated by $(-1,-1,-1), (-1,-1,1)$, hence $(\sigma_2)^\vee_m$ has dimension $1$ and $\cF(U_{\sigma_2})_m=0$.

The cone $\tau^\vee$ is generated by $(0,-1,1), (0,-2,-1), \pm(1,1,0)$. To generate the semigroup $\tau^\vee\cap M$, we need one additional point, for example 
\[ q = \frac{1}{3} ((0,-1,1)+(0,-2,1)) = (0,-1,0).\]
Then $m=2q+(1,1,0)$. The intersection $\tau^\vee \cap (-\tau^\vee + m)\cap M$ is the set
\[ \{ k(1,1,0), q+i(1,1,0), 2q+j(1,1,0)| k,i,j\in\ZZ\}.\]
These points span a plane in $\RR^3$. Since $m$ lies in the interior of $\tau^\vee$, $\Span \tau^\vee_m = \RR^3$.  It follows that $\cF(U_\tau)_m$ has dimension $1$. 

This finishes the proof that $H^1(X,\cF)\neq 0$ when $X= U_{\sigma_1}\cup U_{\sigma_2}$. Next consider $X_\Delta = X\cup Y$, where $Y$ is the toric variety corresponding to the fan $\Delta\setmin\{\sigma_1,\sigma_2, \tau\}$.  Mayer-Vietoris sequence for the open cover is
\[  \ldots \to  H^1(X_\Delta,\cF) \to H^1(X,\cF)\oplus H^1(Y,\cF) \to H^1(X\cap Y,\cF) \to \ldots .\]
We claim that $H^1(X\cap Y,\cF)=0$, which then implies that  $H^1(X_\Delta,\cF) \neq 0$. Indeed, $X\cap Y$ is covered by open affines $U_{\tau_i}$, where $\tau_i$ has dimension $2$. The intersections  $U_{\tau_i} \cap U_{\tau_j}$ for $i\neq j$ are  smooth, hence $\cF$ is zero on these intersections. Using $\check{C}$ech cohomology, we see that  $H^1(X\cap Y,\cF)=0$.
\end{proof}

\subsection{The case of arbitrary complete toric 3-folds.} 
    
Let $X$ be a toric variety corresponding to a $3$-dimensional complete fan $\Delta$. We consider two cases: when every ray in $\Delta$ has at least $4$ neighbours, and when some ray has only $3$ neighbours. We show that in the first case $X$ has a nontrivial line bundle. In the second case $X$ admits a finite surjective toric morphism $X'\to X$, such that $K_0(X')$ is uncountable.

\begin{lemma} Suppose every $1$-dimensional cone $\rho \in\Delta$ is contained in at least $4$ $2$-dimensional cones of $\Delta$. Then there exists a nontrivial integral conewise linear function on the fan $\Delta$.   
\end{lemma}

\begin{proof} This is a simple dimension count. Let $f_i$ be the number of $i$-dimensional cones in $\Delta$. For a cone $\sigma$ of dimension $3$, let $n_\sigma$ be the number of $2$-dimensional faces of $\sigma$, and for $\rho$ of dimension $1$, let $m_\rho$ be the number of $2$-dimensional cones containing $\rho$. 

We represent a conewise linear function on $\Delta$ by its values $a_\rho$ on the marked points $n_\rho$ of rays $\rho$. This gives $f_1$ variables. For every $3$-dimensional cone $\sigma$ there are $n_\sigma-3$ relations. This gives a total of 
\[ \sum_\sigma (n_\sigma-3) = \sum_\sigma n_\sigma  -3 f_3 = 2f_2-3f_3\]
relations. There is a $3$-dimensional space of trivial conewise linear functions, hence if we can prove that 
\[ f_1> 2f_2-3f_3+3,\]
then we can find a rational, or even integral, nontrivial conewise linear function on $\Delta$. 

Using the Euler's identity $f_1-f_2+f_3=2$, we can eliminate $f_3$ and the required inequality becomes
\[ f_2>2 f_1-3.\]
The assumption on the fan is that $m_\rho \geq 4$ for every ray $\rho$, hence
\[ 4 f_1 \leq \sum_\rho m_\rho  = 2f_2.\]
This last inequality implies that $f_2 \geq 2 f_1 > 2f_1-3$.
\end{proof}

\begin{lemma} Let $\Delta$ be a complete fan and $\rho\in\Delta$ a $1$-dimensional cone that is contained in only three $2$-dimensional cones of $\Delta$. Then there exists a sublattice of finite index $N'\subset N$, such that the toric variety $X'$ defined by $(\Delta, N')$ has $K_0(X')$ of uncountable rank.
\end{lemma}

\begin{proof} Let $\tau,\tau_1,\tau_2$ be the $2$-dimensional cones containing $\rho$. By simple geometric considerations one can choose $l\in M$, such that $l\in \Relint \tau^\vee$, but $l\notin \tau_1^\vee$, $l\notin \tau_2^\vee$. We claim that, replacing $N$ by a sublattice $N'\subset N$ and $l$ by a multiple $p\cdot l$ if necessary, we may assume that the cone $\tau$ is nonsingular in $N'$ with primitive generators $v_1,v_2$, and $l$ takes the value $1$ on both $v_1$ and $v_2$. To achieve this, first take the sublattice $N''=N_1\oplus N_2$, where $n\in N\setminus (\Span(\tau)\cap N)$, $N_1 =\Span(n)\cap N$ and $N_2 = \Span(\tau)\cap N$. If $v_1, v_2\in N_2$ are generators of the cone $\tau$ (not necessarily primitive), such that $l(v_1)=l(v_2) = q$, then replace $N_2$ by the sublattice generated by $v_1, v_2$, and further replace $l$ by $l/q$ and let $M''=\Hom(N'',\ZZ)$.
Let now $m=l/2$, let $M' = M''+m \ZZ \subset M''_\RR$ and let $N'\subset N''\subset N$ be the dual lattice. We repeat the computation of $H^1(X',\cF)$ as above, using the lattice point $m \subset M'$, cone $\tau$ and its two $3$-dimensional neighbours $\sigma_1$, $\sigma_2$.
Since $m\notin \sigma_i^\vee$, we get that $\cF(U_{\sigma_i})_m = 0$ for $i=1,2$. It is clear that the points
\[ J:=\tau^\vee\cap (m-\tau^\vee)\cap M'\]
span a rank $2$ lattice $\tilde J\subset M'$. Since $m\in \Relint \tau^\vee$, it follows that $\tilde\Omega^1_X(U_\tau)_m=k.\chi^m\otimes M'$ and thus
\[\cF(U_\tau)_m=k\otimes_\ZZ(M'/\tilde J)\simeq k.\]
The rest of the computation is the same as before.
\end{proof}

\bibliographystyle{plain}
\bibliography{toric-vector}{}

\def\cprime{$'$}
\begin{thebibliography}{1}

\bibitem{CHWW}
G.~Corti{\~n}as, C.~Haesemeyer, M.~Walker, and C.~Weibel.
\newblock The {$K$}-theory of toric varieties.
\newblock {\em Trans. Amer. Math. Soc.}, 361(6):3325--3341, 2009.

\bibitem{Fulton}
William Fulton.
\newblock {\em Introduction to toric varieties}, volume 131 of {\em Annals of
  Mathematics Studies}.
\newblock Princeton University Press, Princeton, NJ, 1993.
\newblock The William H. Roever Lectures in Geometry.

\bibitem{Gubeladze}
Joseph Gubeladze.
\newblock Toric varieties with huge {G}rothendieck group.
\newblock {\em Adv. Math.}, 186(1):117--124, 2004.

\bibitem{Klyachko}
A.~A. Klyachko.
\newblock Equivariant bundles over toric varieties.
\newblock {\em Izv. Akad. Nauk SSSR Ser. Mat.}, 53(5):1001--1039, 1135, 1989.

\bibitem{Payne}
Sam Payne.
\newblock Toric vector bundles, branched covers of fans, and the resolution
  property.
\newblock {\em J. Algebraic Geom.}, 18(1):1--36, 2009.

\bibitem{VV}
Gabriele Vezzosi and Angelo Vistoli.
\newblock Higher algebraic {$K$}-theory for actions of diagonalizable groups.
\newblock {\em Invent. Math.}, 153(1):1--44, 2003.

\end{thebibliography}

\end{document}